\documentclass{tac}
\usepackage{amsmath}

\title {Categories, norms and weights}

\author {Marco Grandis}

\address{Dipartimento di Matematica
\\ Universit\`{a} di Genova
\\Via Dodecaneso 35
\\16146-Genova, Italy }

\eaddress{grandis@dima.unige.it}

\thanks{Work supported by MIUR Research Projects.}



\date{}

\def \c {\colon}

\def \q {\qquad}
\def \qq {\qquad \qquad}
\def \ndt {\noindent}
\def \bu {{\scriptscriptstyle\bullet}}

\def \ti {\!\times\!}
\def \te {\otimes}
\def \jn {\vee}
\def \mt {\wedge}
\def \ci {\!\mathbin{\raise.3ex\hbox{$\scriptscriptstyle\circ$}}\!}
\def \iso {\: \cong \:}
\def \eq {\! \sim \!}
\def \iff {\: \Leftrightarrow \:}
\def \for {\; \mbox{ for } }

\def \min {{\rm min}}
\def \sup {{\rm sup}}
\def \inf {{\rm inf}}
\def \ln {{\rm ln}}
\def \exp {{\rm exp}}
\def \hom {{\rm hom}}
\def \Hom {{\rm Hom}}
\def \Re {{\rm Re}}
\def \Im {{\rm Im}}
\def \arg {{\rm arg}}

\def \op{^{{\rm op}}}

\def \la {\lambda}
\def \ph {\varphi}

\def \Set {{\bf Set}}
\def \wSet {{\rm w} {\bf Set}}
\def \wpSet {{\rm w}^+ {\bf Set}}
\def \wbSet {{\rm w}^{\bu} {\bf Set}}
\def \wpMon {{\rm w}^+ {\bf Mon}}
\def \wbMon {{\rm w}^{\bu} {\bf Mon}}
\def \wAb {{\rm w} {\bf Ab}}
\def \wpAb {{\rm w}^+\! {\bf Ab}}
\def \dAb {{\rm d} {\bf Ab}}
\def \wRng {{\rm w} {\bf Rng}}

\def \wpCat {{\rm w}^+ {\bf Cat}}
\def \wbCat {{\rm w}^{\bu} {\bf Cat}}

\def \bt {{\bf 2}}
\def \W {{\bf w}}
\def \v {{\bf v}}
\def \p {{\bf p}}
\def \r {{\bf r}}
\def \Wp {{\bf w}^+}
\def \Wb {{\bf w}^{\bu}}

\def \sing  {\{*\}}
\def \N {{\bf N}}

\def \R {{\bf R}}
\def \C {{\bf C}}
\def \I {{\bf I}}

\def \iy {\infty}
\def \zety {[0, \infty]}

\def \wp {w$^+$}
\def \wb {w$^{\bu}$}
\def \wZ {{\rm w}{\bf Z}}
\def \wQ {{\rm w}{\bf Q}}
\def \wR {{\rm w}{\bf R}}

\begin{document}

\maketitle

\begin{abstract}

	The well-known Lawvere category $ \zety $ of extended real positive numbers comes with a monoidal closed structure where the tensor product is the sum. But $ \zety $ has another such structure, given by multiplication, which is *-autonomous and a CL-algebra (linked with classical linear logic).

	{\it Normed sets}, with a norm in $ \zety$, inherit thus two symmetric monoidal closed structures, and categories enriched on one of them have a `subadditive' or `submultiplicative' norm, respectively. Typically, the first case occurs when the norm expresses a cost, the second with Lipschitz norms.
	
	This paper is a preparation for a sequel, devoted to {\it weighted} algebraic topology, an enrichment of {\it directed} algebraic topology. The structure of $ \zety$, and its extension to the complex projective line, might be a first step in abstracting a notion of {\it algebra of weights}, linked with physical measures.

\end{abstract}

\section*{Introduction}\label{Intro}

	A category can be equipped with a {\it (sub)additive} norm satisfying, for all objects $ X $ and all pairs of composable arrows $ f, g $
%
    \begin{equation}
|1_X| =  0,   \qq   |gf|  \le  |f| + |g|,
    \label{0.1}  \end{equation}
or also with a {\it (sub)multiplicative} norm, satisfying:
%
    \begin{equation}
  |1_X|  \le  1,   \qq   |gf|  \le  |f|.|g|.
    \label{0.2}  \end{equation}
    
    	The first case appears when the norm expresses a `cost' (length, duration, price, energy,...) which can `at worst' be added in a composition (typically, in a concatenation of paths). The second is usual with Lipschitz norms, where the norm expresses a scale factor (or, rather, a best bound for that), which can `at worst' be multiplied in a composition. A normed additive category (like normed vector spaces and bounded linear mappings) makes use of both aspects: its hom-sets are abelian groups equipped with an {\it additive} norm, but composition is {\it multiplicative}. The `same' happens in a {\it normed ring} - a normed additive category on one object, or a monoid object in the category of additively weighted abelian groups, with multiplicative tensor product. The purpose of this article is to investigate such aspects and fix a coherent terminology.
	
	As shown in the well-known article of Lawvere on (generalised) metric spaces \cite{Lw}, the additive notion (1) is based on the category of extended positive real numbers, $  \la  \in \zety$, with arrows $  \la \ge  \mu$, equipped with a strict symmetric monoidal closed structure, which we write $ \Wp$: the tensor product is the sum $  \la +  \mu $ and the internal hom is defined by truncated difference, $ \hom^+( \mu,  \nu) = 0  \jn ( \nu -  \mu)$. One derives from $ \Wp $ the symmetric monoidal closed category $  \wpSet $ of {\it normed sets} (\ref{1.4}), written $ {\cal S}(\R) $ in \cite{Lw} and more explicitly described in a paper by Betti-Galuzzi \cite{BG}. A category enriched on the latter (cf. \ref{3.1}) has an additive norm, as in (\ref{0.1}). 

	But the same category $ (\zety, \ge) $ can be equipped with a {\it multiplicative} tensor product 
$  \la . \mu$.
 {\it Provided we define} $ 0.\iy = \iy $ (so that tensoring by 0 preserves the initial object $ \iy$), the latter is again a strict symmetric monoidal closed structure $ \Wb$: the internal hom is 
$ \hom^ {\bu}( \mu,  \nu) =  \nu/ \mu$, 
where the `undetermined forms' $ 0/0 $ and $ \iy / \iy $ are defined to be 0. (This choice comes from privileging the direction $  \la \ge  \mu$, which is necessary if we want to view (2) as an expression of enrichment; cf. \ref{1.2}.) Multiplication gives a {\it multiplicative} symmetric monoidal closed structure $  \wbSet$, on the category of normed sets; enrichment on the latter means a multiplicative norm, cf. \ref{3.2}. The new multiplicative structure $ \Wb $ is *-autonomous (Barr \cite{Ba}), with involution $ 1/  \mu $ (\ref{1.2}). (Lipschitz norms are viewed in \cite{Lw} in a different way, based on endofunctors of $ \Wp$, cf. \ref{3.5}.)
	
	Now, these `norms' are `generalised norms', in the same way as Lawvere metrics are generalised ones:

\ndt  (a) they can take an infinite value,

\ndt  (b) they are not `Hausdorff': $ |a| = 0 $ does not imply $ a = 0 $ (when this makes sense),

\ndt  (c) they are `directed': we do not assume $ |a| = |- a| $ (when this makes sense). 
		
	This is why we prefer to speak of {\it weights} (or {\it costs}) rather than norms: a {\it weighted abelian group}, as defined here (\ref{2.2}), is a much weaker notion than a normed abelian group in the classical sense. As a typical example, the free weighted abelian group on one element is the group of integers $ \wZ $, equipped with the weight where $ |k| = k $ for $ k \ge 0 $ and otherwise $ |k| = \iy $ (\ref{2.2}); the weighted field $ \wR $ of reals has a similar weight (\ref{2.5}). These examples also exhibit how {\it weights (can) have a directed character}: a weighted abelian group has an associated preorder, where the positive cone is given by the {\it attainable} elements (of finite weight; cf. \ref{2.4}).
	
	This paper is a preparation for a sequel where, starting from some works in `directed algebraic topology' (cf. \cite{G2, G3} and references therein), we shall develop a study of `weighted algebraic topology', where `weighted spaces' (for instance, the generalised metric spaces of Lawvere) can be studied - via homotopy or homology - with {\it additively weighted} categories (\ref{3.1}) or with {\it additively weighted} abelian groups and their {\it multiplicative} tensor product. A first study of weighted homology for weighted cubical sets can be found in \cite{G4}, showing links with noncommutative geometry.
	
	Finally, it may be interesting to pay some reflection to the {\it algebraic} structure which makes $ \zety $ an `algebra of measures' or `weights'. As briefly considered in \ref{1.8}, $ \zety $ is a commutative semiring, equipped with an involution which transforms sum and multiplication into other operations. Such a structure also exists in the complex projective line $ P^1\C$, whose additive structure and involution $ z^{-1}$ yield the calculus of {\it impedances} and {\it admittances}, in RLC networks. But $ \zety $ has an {\it order}, which allows us to treat it as a category (cf. \ref{1.2}), and it is not clear if and how this aspect should be partially transferred to $ P^1\C $ (cf. \ref{1.8}).

\ndt  {\it Acknowledgements}. The author is indebted with F.W. Lawvere, R. Betti, G. Rosolini for various helpful discussions.

\section{Weights and weighted sets}\label{1}

	Weighted sets, equipped with a weight function $ X \to \zety$, have two tensor products, derived from two tensor products on the category of weights $ \zety$, i.e. addition and multiplication. All these structures are symmetric monoidal closed. We end with some speculation on `algebras of weights' (\ref{1.8}).

\subsection{Real weights}\label{1.1}
	Quantities are generally measured by positive real numbers. It is convenient to use the interval $ \zety $ of {\it extended} positive real numbers, including $\iy$, which will be called {\it real weights}, or just {\it weights}. For instance, if we are considering the resistance $ R $ of a conductor, the perfect conductor $ (R = 0) $ should have equal rights to be included as the perfect insulator $ (R = \iy)$; also because the conductance $ G = R^{-1} $ reverses such values.

	The interval $ \zety $ has a rich structure, of relevance for measures. First, it is totally ordered, which allows us to compare quantities of the same species, and a complete lattice (thanks to including $\iy$), obviously distributive.
	
	Second, $ \zety $ has a sum (which measures a disjoint union of quantities of the same species) and a multiplication (corresponding to the product of physical quantities). It is a commutative {\it semiring}, letting
%
    \begin{equation}
\la + \iy  =   \la.\iy  =   \la \jn \iy  =  \iy   \qq\qq   ( \la  \in \zety).
    \label{1.1.1}  \end{equation}

	Since we want to use both operations, and some physical intuition here and there, we should insist on the fact that these weights are `pure quantities': they do not stand for physical quantities but for their measures, with respect to a fixed system of units. Let us also note that $ \iy $ acts in the obvious way, except for the choice $ 0.\iy = \iy $ of the `undetermined form', which is motivated below (\ref{1.2}). As to terminology, a {\it semiring} is here an (additive) abelian monoid equipped with a multiplication which is associative, has a unit and distributes over addition; note that we do not (cannot) require that 0 be absorbent for multiplication: the present notion is slightly weaker than a {\it rig} (which is a monoid in the monoidal category of abelian monoids, cf. \cite{Sc}). The prime example is the semiring $ \N $ of natural numbers, which is the free such algebra on the empty set. 

	Finally, there is in $ \zety $ a (well defined) involution $  \la \mapsto  \la^{-1}$, whose physical interest is also obvious. This involution transforms the sum into another operation, which we call {\it harmonic sum} 
%
    \begin{equation}
 \la *  \mu  =  ( \la^{-1} +  \mu^{-1})^{-1},
    \label{1.1.2}  \end{equation}
of common use for $l_p$-norms, but also in geometrical optics and for resistor networks: in the latter case, composition in series of resistors leads to the sum of resistances, while composition in parallel leads to the sum of conductances, and the harmonic sum of resistances. (This interpretation suggests an extension to complex numbers, measuring impedances and admittances, see \ref{1.8}.)

	Similarly, the involution transforms the multiplication into a commutative monoid operation
%
    \begin{equation}
\la  \bu  \mu  =  ( \la^{-1}. \mu^{-1})^{-1},
    \label{1.1.3}  \end{equation}
which coincides with multiplication, except that it gives `the other choice' for the undetermined form: $ 0 \bu\iy = 0$.

	One can also note that all translations $  \la + ( \; ) \c \zety \to \zety $ preserve arbitrary meets (while, generally, they do not preserve the empty join, 0); thus, reversing the order, $ (\zety, \ge, +) $ is a commutative unital quantale. The same holds with the multiplicative structure $ (\zety, \ge, .)$, with definition (\ref{1.1.1}). On the other hand, the involution $  \la \mapsto  \la^{-1} $ does not make any of these structures into an `involutive quantale', according to the notion in use (cf. Mulvey-Pelletier \cite{MP}). The categorical counterpart of these facts is the closure of the two monoidal structures, $  \la +  \mu $ and $  \la. \mu$, considered below (\ref{1.2}).

\subsection{The category of weights}\label{1.2}
	Our weight functions - for sets, abelian groups or categories - will take values in the commutative ordered semiring $ \zety $ considered above. But we need to make it into a category.
	
	As in Lawvere \cite{Lw}, we use the {\it opposite} category $ \W $, with morphisms $  \la \ge  \mu $ (which is necessary to treat normed categories as enriched ones). This category has all limits and colimits (with trivial equalisers and coequalisers, of course)
%
    \begin{equation}
    \begin{array}{lllll}
\mbox{ product:}  & \sup  \la_i,
&\q&   \mbox{ terminal object:} &  0,
\\
\mbox{ sum:} &  \inf \la_i,
&\q&  \mbox{  initial object:}  &  \iy.
    \label{1.2.1}\end{array}
    \end{equation}

	We will write $ \Wp = (\W, +, 0) $ this category equipped with the strict symmetric monoidal closed structure defined by the {\it additive} tensor product, $  \la +  \mu$. The internal hom is given by {\it truncated subtraction}, and will be written as a difference (as in \cite{Lw}): 
%
    \begin{equation}
\la +  \mu \ge  \nu  \;   \iff  \;    \la \ge \hom^+( \mu,  \nu)  =   \nu -  \mu,
\qq	( \nu -  \mu  =  0  \jn  ( \nu -  \mu)).
    \label{1.2.2}  \end{equation}

	Note that an `undetermined form' appears, and gets a precise value: $ \iy - \iy = 0 $ (since $  \la + \iy \ge \iy $ is always true, as $  \la \ge  0$).
	
	We will write $ \Wb = (\W, ., 1) $ the same category, equipped with the strict symmetric monoidal closed structure defined by the {\it multiplicative} tensor product $  \la. \mu$; recall that we have chosen $ 0.\iy = \iy$, in \ref{1.1}. (Thus, tensoring with {\it any} element $  \la $ preserves the initial object $ \iy$, for the `direction' $  \la \ge  \mu$. But of course the opposite category $ \W \op = (\zety, \le) $ is monoidal closed with tensor product $  \la \bu \mu$, and $ \iy \bu 0 = 0$.)
	
	Now, the internal hom for $ \Wb $ is given by division, with `undetermined forms' $ 0 / 0 = \iy / \iy = 0 $ (as required by the adjunction, after the previous choice)
%
    \begin{equation}
    \begin{array}{l}
\la. \mu \ge  \nu    \;   \iff  \;    \la \ge  \hom^ \bu ( \mu,  \nu)  =   \nu /  \mu,
\\	( \nu / \iy = 0   \; \mbox{  for all  }   \nu;   \q    \nu / 0 = \iy   \for    \nu > 0,   \q   0 / 0  =  0).
    \label{1.2.3}\end{array}
    \end{equation}

	This is a *-autonomous category \cite{Ba}, with dualising object the multiplicative identity 1 and involution $  \la \mapsto \hom^ \bu ( \la, 1) =  \la^{-1}$. (It is thus a CL-algebra, i.e. an algebra for classical linear logic according to Troelstra's book \cite{Tr}, or a commutative cyclic Grishin algebra according to Lambek \cite{Lm}.)
	
	We have already remarked that the derived operation $  \la \bu  \mu = ( \la^{-1}. \mu^{-1})^{-1} $ `nearly' coincides with the multiplication (\ref{1.1}): the *-autonomous structure is `nearly' compact. One can also note that all finite $ \xi > 0 $ are dualising elements, i.e. provide an involution $  \la \mapsto \hom^ \bu ( \la, \xi) = \xi /  \la$; but $ \xi = 1 $ is the unique choice giving a `nearly compact' structure.
	
	Note also that the values of the `undetermined forms', $ 0 / 0 = 0 = \iy / \iy$, agree with the following identity (holding in every *-autonomous category, in the right-hand form)
%
    \begin{equation}
	 \la /  \mu  =  ( \mu. \la^{-1})^{-1}	   \qq   (\hom( \mu,  \la) \iso ( \mu\te \la^*)^*).
    \label{1.2.4}    \end{equation}

	The two tensor products in the category $ \W $ have an interesting interplay, already examined above from the algebraic point of view. Multiplication distributes on addition, and the involution $  \la \mapsto  \la^{-1} $ of the *-autonomous multiplicative structure transforms $ \Wp $ into an anti-isomorphic symmetric monoidal closed category $ \W ^*$: the opposite category, equipped with the harmonic sum $  \la *  \mu = ( \la^{-1} +  \mu^{-1})^{-1} $ (cf. (\ref{1.1.2})). These structures will be further examined in \ref{1.8}.

\subsection{Truth values}\label{1.3}
	Let $ \v = ( \{0, \iy \}, \ge) $ denote the full subcategory of $ \W $ on the objects $ 0, \iy$. Note that, in this subcategory, the cartesian product $  \la \jn  \mu $ coincides with both tensor products considered above, $  \la +  \mu $ and $  \la. \mu $, in (\ref{1.1.1}).
	
	The category $ \v $ (`v' for verity) is isomorphic to the boolean algebra $ \bt = ( \{0, 1 \}, \le) $ of {\it truth-values}, a cartesian closed category with cartesian product $ p \mt q = p.q$. The covariant embedding (contravariant with respect to the natural orders)
%
    \begin{equation}
M \c \bt \to \W,   \qq   M(0)  =  \iy,   \q   M(1)  =  0,
    \label{1.3.1}    \end{equation}
transforms the (cartesian) product in $ \bt $ into the three tensor products of $ \W $ (which coincide in $ v$). Moreover, $ M $ has left and right adjoint
%
    \begin{equation}
P  \dashv   M  \dashv   Q,   \qq   P( \la) = 1  \iff   \la < \iy,   \qq   Q( \la) = 1  \iff   \la = 0.
    \label{1.3.2}    \end{equation}

\subsection{Weighted sets}\label{1.4}
A {\it weighted set}, or {\it w-set}, will be a set $ X $ equipped with a {\it weight}, or {\it cost function}, consisting of an arbitrary mapping
%
    \begin{equation}
w_X \c X \to \zety,
    \label{1.4.1}    \end{equation}
also written $ w$, or $ | - |_X$, or $ | - |$. We shall say that an element of $ X $ is {\it free}, {\it attainable} or {\it unattainable} when, respectively, its cost is 0, finite or $\iy$.

	A (weak) {\it contraction} $ f \c X \to Y$, or {\it w-map}, or map of w-sets, has $ |f(x)| \le |x|$, for all $ x  \in X$. $ \wSet $ will denote the category of these {\it weighted sets and contractions}; an isomorphism is thus a bijective isometry: $ |f(x)| = |x|$, for all $ x$. This category has all limits and colimits, constructed as in $ \Set $ and equipped with a suitable weight (strictly determined).

	Thus, a product $ \prod X_i $ and a sum $ \sum X_i $ (where $ X_i $ has weight $ | - |_i$) have the following weights  
%
    \begin{equation}
    \begin{array}{lll}
|(x_i)|  =  \sup_i |x_i|_i   &\qq&   ((x_i)  \in \prod X_i),
\\     |(x, i)|  =  |x|_i  &\qq&   (x  \in X_i),
    \label{1.4.2}    \end{array}
     \end{equation}
while a {\it weighted subset} has the restricted weight, and a {\it quotient} $ X/\eq \;$ has the induced one
%
    \begin{equation}
	|\xi|  =  \inf \{|x| \;  | \;  x \in \xi \}   \qq   (\xi \in X/\eq).
    \label{1.4.3}    \end{equation}

	Plainly, infinite products exist because we are allowing an infinite weight. By the same reason, the forgetful functor $ B_\iy \c \wSet \to \Set $ has a left adjoint $ w_\iy S$, which equips the set $ S $ with the {\it discrete weight}, always $ \iy$, and a right adjoint $ w_0S$, with the {\it codiscrete weight,} always zero.
	
	More generally, the functor $ w _\la \c \Set \to \wSet $ which equips a set with the constant weight $  \la \in \zety $ has for right adjoint the (representable) $ \la$-{\it ball} functor
%
    \begin{equation}
B _\la \c \wSet \to \Set,   \qq   B _\la(X)  =   \{x \in X \;  | \;   |x| \le  \la \}  =  \wSet(w _\la \sing, X),
    \label{1.4.4}    \end{equation}
and there is a chain of adjunctions: $ w_\iy \dashv B_\iy \dashv w_0 \dashv B_0$.

\subsection{The additive tensor product}\label{1.5}
	This structure, hinted at in \cite{Lw}, is explicitly described in \cite{BG}.
	
	We shall write $  \wpSet $ the closed symmetric monoidal category of w-sets, equipped with the {\it additive} tensor product $ X \te_0 Y$, derived from the tensor product of $ \Wp$. It is given by the cartesian product $ |X| \ti |Y| $ of the underlying sets, with the following {\it additive} weight on a pair $ x \te_0 y $ (written thus to avoid confusion with the cartesian product, where $ |(x, y)| = |x| \jn |y|$)
%
    \begin{equation}
|x \te_0 y|  =  |x| + |y|.
    \label{1.5.1}    \end{equation}

	The identity of the tensor product is the 0-weighted singleton $ w_0 \sing$, and the representable functor produced by the latter is the zero-ball functor $ B_0 $ (\ref{1.4.4})
%
    \begin{equation}
B_0 \c \wSet \to \Set,   \qq   B_0(X)  =   \{x \in X \; | \; |x| = 0 \}  =  \wSet(w_0 \sing, X).
    \label{1.5.2}    \end{equation}

	The internal hom is the set of {\it all} mappings, equipped with the {\it additive weight}, or {\it truncated-difference weight}
%
    \begin{equation}
    \begin{array}{l}
\Hom_0 (Y, Z)  =  \Set(|Y|, |Z|),
\\[2pt]  
  |h|_0  =  \sup_y (|h(y)| - |y|)  =
  \min  \{ \la \in \zety \; | \; \mbox{ for all }  y \in Y, \;  |h(y)| \le  \la + |y| \},
\\[2pt]  
 |h|_0 \le  \la  \;  \iff  \mbox{  for all }  y \in Y, \;  |h(y)| \le  \la + |y|.
    \label{1.5.3}    \end{array}
     \end{equation}

	Thus, the usual bijection $ \Set(X \ti Y, Z) = \Set(X, \Set(Y, Z)) $ which identifies $ f \c X \ti Y \to Z $ with $ g \c X \to \Set(Y, Z) $ under the condition $ f(x, y) = g(x)(y)$, provides two isometries
%
    \begin{equation}
    \begin{array}{lll}
\Hom_0(X \te_0 Y, Z)  =  \Hom_0(X, \Hom_0(Y, Z)),   &\q&   |f|_0  =  |g|_0,
\\[2pt]  
  \wSet(X \te Y, Z)  =  \wSet(X, \Hom_0(Y, Z)).
    \label{1.5.4}    \end{array}
     \end{equation}

	And of course, the zero-ball functor $ B_0, $ applied to the internal hom, gives back the w-maps
%
    \begin{equation}
B_0(\Hom_0(Y, Z))  =  \wSet(Y, Z).
    \label{1.5.5}    \end{equation}

\subsection{The multiplicative tensor product}\label{1.6}
	Similarly, $ \wSet $ has a second important closed symmetric monoidal structure, which we shall write $  \wbSet$. The {\it multiplicative} tensor product $ X \te_1 Y $ is given again by the cartesian product $ |X| \ti |Y| $ of the underlying sets, with the {\it multiplicative} weight derived from the tensor product of $ \Wb $ (recall that $ 0.\iy = \iy$)
%
    \begin{equation}
|x \te_1 y|  =  |x|.|y|.
    \label{1.6.1}    \end{equation}

	The identity is the 1-weighted singleton $ w_1 \sing$, and its representable functor is the unit-ball:
%
    \begin{equation}
B_1 \c \wSet \to \Set,	   \qq   B_1(X)  =   \{x \in X | |x| \le 1 \}  =  \wSet(w_1 \sing, X).
    \label{1.6.2}    \end{equation}

	The internal hom is the set of all mappings equipped with the {\it multiplicative weight}, or 
{\it Lipschitz weight}, i.e. the least Lipschitz constant of a mapping (possibly $\iy$)
%
    \begin{equation}
    \begin{array}{l}
\Hom_1(Y, Z)  =  \Set(|Y|, |Z|),
\\[2pt]  
|h|_1  =  \sup_y (|h(y)| / |y|)  =
  \min  \{ \la \in \zety \; | \; \mbox{ for all }  y \in Y, \;  |h(y)| \le  \la.|y| \},
\\[2pt]  
 |h|_1  \le  \la  \;  \iff  \mbox{  for all }  y \in Y, \;  |h(y)| \le  \la.|y|.
    \label{1.6.3}    \end{array}
     \end{equation}
	Again, the exponential law in $ \Set $ provides two isometries (notation as in \ref{1.5.4})
%
    \begin{equation}
    \begin{array}{lll}
\Hom_1(X \te_1 Y, Z)  =  \Hom_1(X, \Hom_1(Y, Z)),   &\q&   |f|_1  =  |g|_1,
\\[2pt]  
  \wSet(X \te_1 Y, Z)  =  \wSet(X, \Hom_1(Y, Z)).
    \label{1.6.4}    \end{array}
     \end{equation}
and the unit-ball functor $ B_1$, applied to the internal hom, gives back the w-maps
%
    \begin{equation}
B_1(\Hom_1(Y, Z))  =  \wSet(Y, Z).
    \label{1.6.5}    \end{equation}

\subsection{Probabilistic and relative weights}\label{1.7}
The additive structure $ \Wp $ is isomorphic to the category $ \p = [0, 1] $ of {\it probabilistic weights}, with morphisms $ p \le q$, via
%
    \begin{equation}
\la  =  - \, \ln(p),   \qq   p  =  \exp(-  \la)  
    \label{1.7.1}    \end{equation}
(and {\it anti}-isomorphic as a lattice, with respect to the natural orders). The category $ \p $ has thus an isomorphic structure, with internal hom by truncated division:
%
    \begin{equation}
    \begin{array}{llll}
\mbox{product and sum:}  \q &  \inf (p_i),   & \sup (p_i),  \q
\\
\mbox{tensor product:}  &  p.q,
\\
\mbox{internal hom:} \q  &  \hom(q, r) = 1  \mt  r/q,
\\
\mbox{adjunction:}   &   pq \le r   \iff   p \le 1  \mt  r/q. \q
    \label{1.7.2}\end{array}
    \end{equation}
	Moreover, $ p $ contains the category $ \bt$, on which cartesian and tensor product coincide.
	
	On the other hand, the multiplicative structure $ \Wb $ is isomorphic to the category $ \r = [- \iy, \iy] $ of {\it relative weights}, with morphisms $ x \ge y$, via
%
    \begin{equation}
x  =   \ln( \la),   \qq   p  =  \exp(x). 
    \label{1.7.3}    \end{equation}

	The category $ \r $ has thus an isomorphic structure of *-autonomous category, with
%
    \begin{equation}
    \begin{array}{lll}
\mbox{product and sum:}  &  \sup (x_i), & \inf (x_i),
\\
\mbox{tensor product:}  \q &  x+y,  &    (- \iy + \iy  =  \iy),
\\
\mbox{internal hom and involution:}  \;  & \hom(y, z)  =  z - y, \q   &   - x
    \label{1.7.4}\end{array}
    \end{equation}

\subsection{From cubical monoids to algebras of weights}\label{1.8}
	Before going on with the main goals, let us pay some attention to the structure of $ \zety$, which makes it an `algebra of weights', adequate to express measures of physical quantities. We shall distinguish some properties of this kind, without giving a precise definition of an algebra of weights, which would require a deeper study.

	Let us first note that measures need not be confined to the real line; complex numbers are also used, e.g. in the analysis of electric networks (or of their mechanical equivalents). The {\it purely algebraic} structure examined above, with main operations 
$  \la +  \mu, \,  \la. \mu,  \, \la^{-1} $ 
extends, with the same properties, to the complex projective line $ P^1\C = \C \cup  \{\iy \}$ - and more generally to the projective line $ P^1F $ on any commutative field $ F $ (always defining $  \la + \iy = \iy =  \la.\iy$, for all $  \la$).
	
	Extending the previous interpretation for resistor networks, the additive structure of $ P^1\C$, together with the involution $ z^{-1}$, formalises the calculus of {\it impedances} and {\it admittances}, for networks of resistors, inductors and capacitors in steady sinusoidal state (cf. \cite{Vv}). In this situation, the complex number $ Z = R + iX $ represents the impedance of an RLC network, for an alternate current of (fixed) frequency $ f $ and angular speed $ \omega = 2pf$; the real part is the resistance $ R \ge 0 $ while $ X $ is the reactance of our device
%
    \begin{equation}
X  =  \omega L - (\omega C)^{-1}
    \label{1.8.1}    \end{equation}
which results of its inductance $ L $ and capacitance $ C$. Sum and harmonic sum of impe-dances still agree with composition in series and parallel, respectively.

	For this interpretation, it makes sense to restrict $ P^1\C $ to $ \C^+$, the complex numbers with real part $ \ge 0 $ (including $\iy$), which are closed under the structure we are considering - sum, involution and harmonic sum - but not under product. However, the multiplicative structure of (the whole) $ P^1\C $ is also of interest for physical measures, e.g. in the analysis of electrical networks in alternate current.

	The goal of abstracting a notion of `algebra of measures', or `weights', might begin with the following steps.
	
\ndt (a) Let us start from a {\it dioid}, or {\it cubical monoid}, a structure of interest for homotopy and standard intervals, introduced in \cite{G1}: it is a set equipped with two structures of monoid, such that the unit of each operation is an absorbent element for the other. Typically, a lattice (with minimum and maximum) has such a structure, but a cubical monoid need not satisfy the idempotence laws, nor the general absorption laws. (The name comes from links with cubical sets. A cubical monoid in a category of endofunctors is called a cubical monad \cite{G1}, typically the structure of a cylinder functor; thus, an augmented simplicial set is to a {\it cubical set} what a monad is to a {\it cubical monad}, what a monoid is to a {\it cubical monoid}.)

\ndt (b) An {\it involutive cubical monoid} \cite{G1} has an involution turning each structure into the other. Equivalently, one can give a monoid $ (A, +, 0) $ equipped with an involution $ (-)^* $ such that the element $ \iy = 0^* $ is absorbent for the sum; the second operation is then defined as $ x * y = (x^* + y^*)^*$.

	The standard interval $ \I = [0, 1] $ has two such structures of interest for homotopy: the usual one, as a (totally ordered) lattice with involution $ t^* = 1 - t $, and a non-idempotent structure $ (\I, ., 1, ^*) $ with ordinary multiplication and the same involution, which is of interest for {\it smooth} homotopies in differentiable manifolds (while the lattice operations are not smooth). The {\it ordered interval} has the same structures without the involution. An MV-algebra, used in `multi-valued logic' \cite{Ch, CDM} is a commutative involutive cubical monoid satisfying one axiom more: $ (x^* + y)^* + y = (y^* + x)^* + x$.

\ndt (c) Now, let us define a {\it cubical semiring} as a set equipped with two structures of semiring, such that the zero of each structure is an absorbent element for both operations of the other structure, while the multiplicative units coincide.

\ndt (d) Further, an {\it involutive cubical semiring }will also have an involution turning each structure into the other. Again, one can equivalently give a semiring $ (A, +, 0, ., 1) $ with an involution $ ( - )^* $ leaving 1 fixed and such that the constant $ \iy = 0^* $ is absorbent for sum and product
%
    \begin{equation}
x^{**}  =  x,	\q	0^*  =  \iy,    \q  1^*  = 1,
\q   x + \iy  =  x.\iy  =  \iy.x  =  \iy.    
   \label{1.8.2} \end{equation}
	The dual operations 
$ x * y = (x^* + y^*)^*,  \;   x \bu y = (x^*.y^*)^* $
form then an anti-isomorphic structure. In this sense, $ \zety $ and every projective line are involutive cubical semirings, while $ \C^+ $ is a sub-{\it involutive cubical monoid} of $ P^1\C$.
	
	Now, $ \zety $ is a {\it totally ordered} commutative cubical semiring: the semiring structure agrees with its natural order, while the involution reverses it. Something can be said about order properties of $ P^1\C$, provided we separate the additive structure from the multiplicative one, with two different (partial) orders, both extending the total order of $ \zety$.

\ndt (e)  First, one can consider the category with objects in $ \C^+ $ (complex numbers with real part $\ge 0$, including $\iy$) and morphisms $ z \ge^+ w$, meaning that $ \Re(z) \ge \Re(w) $ and $ \Im(z) = \Im(w) $ (with $ \iy $ as an initial object, i.e. a maximum for $ \le^+$). Extending the additive structure of $ \zety$, this category is strict symmetric monoidal closed, with additive tensor product $ z + w $ and $ \hom^+(w, w') = w' - w$, by {\it truncated difference on the real part}. The (partial) order we are considering is - presumably - consistent with the interpretation of our weights as impedances $ R + iX$, where only the real part $ R $ `dissipates' energy. This structure should likely be enriched with the (contravariant) involution $ z^{-1}$, which turns sum into harmonic sum.

\ndt (f)  Second, we have the category with objects in $ P^1\C $ and morphisms $ z \ge^ \bu  w$, meaning that $ |z| \ge |w| $ and $ \arg(z) = \arg(w) $ (with $ \iy $ as an initial object and 0 as a terminal). Extending now the multiplicative structure of $ \zety$, this category is a *-autonomous, with tensor product 
$ z.w$  and  $\hom^ \bu (w, w') = w'/w $ (undetermined forms as in the real case) and dualising object 1.

\section{Weighted algebraic structures}\label{2}

	We fix the terminology for weighted algebraic structures, to be used in the sequel.

\subsection{Weighted monoids}\label{2.1}
An {\it additively weighted monoid} $ A$, or {\it \wp-monoid}, will be a monoid object in the monoidal category $  \wpSet$. Thus, it is a monoid (in additive notation) equipped with a {\it weight function}, written $ w $ or $ | - |$, taking values in $ \zety$ and such that
%
    \begin{equation}
|0|  =  0,   \qq   |a + b|  \, \le \,   |a| + |b|.
   \label{2.1.1} \end{equation}

	In the category $ \wpMon $ of such objects, a {\it morphism} is a contracting homomorphism: $ |f(a)| \le |a|$, for all elements $ a $ of the domain.

	Similarly, a {\it multiplicatively weighted monoid} $ A$, or {\it \wb-monoid}, is a monoid in $  \wbSet$. Writing the operation as a product, the weight satisfies now the following axioms
%
    \begin{equation}
|1|  \le  1,   \qq   |a.b|  \le  |a|.|b|,
   \label{2.1.2} \end{equation}
and we have a category $ \wbMon$, with contracting homomorphisms.

\subsection{Weighted abelian groups}\label{2.2}
	Weighted directed homology, to be studied in a sequel, will take values in the category $ \wAb $ of {\it weighted abelian groups}, i.e. the full subcategory of $ \wpMon $ formed of the additively weighted monoids which are abelian groups. We shall not use the category of {\it multiplicatively} weighted abelian groups, which is why we do not insist in writing $ \wpAb$ .
	
	Note that, for $ n \in \N$, we only have $ |n.a| \le n.|a|$. Note also that we do not require $ |- a| = |a|$; we want a directed notion, able to distinguish `negative elements' by means of an infinite cost (cf. \ref{2.3}). For instance, $ \wZ $ will be the group of integers with
%
    \begin{equation}
w(k)  =  k,   \mbox { for }  k \ge 0,    \qq    w(k)  =  \iy,    \mbox { for }  k < 0.
   \label{2.2.1} \end{equation}

	In any weighted abelian group, we have
%
    \begin{equation}
|ka|  \le  w(k).|a|,    \qq    | \sum  k_i a_i |  \;   \le \;    \sum   w(k_i).|a_i|.
   \label{2.2.2} \end{equation}

	The category $ \wAb $ has all limits and colimits, computed as in $ Ab $ and equipped with a suitable weight (as for w-sets). The tensor product $ A \te B $ of $ Ab $ can be lifted to $ \wAb$, with a multiplicative weight on an element  $ \xi \in A \te B $
%
    \begin{equation}
|\xi|  =  \inf \,  \{\,  \sum  |a_i|.|b_i|  \; \;   |    \;\;    \xi = \sum  a_i\te b_i \, \}    \qq  |a \te  b|  \le  |a|.|b|.
    \label{2.2.3}\end{equation}

	In fact, to prove that $ | \xi + \eta | \le | \xi| + |\eta |$, take any pair of expressions of $ \xi $ and $ \eta$
%
    \begin{equation}
\xi  \; =  \,  \sum_{i=1}^n a_i\te b_i,    \qq    \eta \; =  \sum _{i=n+1}^m  a_i\te b_),
    \label{2.2.4}\end{equation}
and note that
%
    \begin{equation}
| \xi + \eta | \; = \;     |    \sum_{i=1}^m a_i\te b_i     |   
\;  \le \;     \sum _{i=1}^m  |a_i|.|b_i|  \;  = \; 
  \sum_{i=1}^n  |a_i|.|b_i| \, + \sum _{i=n+1}^m  |a_i|.|b_i|.
    \label{2.2.5}\end{equation}

	This object $ A \te B $ solves the universal problem for bi-homomorphisms $ \ph \c A \ti B \to C $ such that $ |\ph(a, b)| \le |a|.|b|$. It produces a closed symmetric monoidal structure: the internal hom $ \Hom(B, C) $ is the abelian group of {\it all} homomorphisms of the underlying abelian groups, with the Lipschitz weight (\ref{1.6.3}).
	
	The unit of the tensor product is $ \wZ$, as defined above (and essentially proved in (\ref{2.2.2}). The representable functor $ \wAb(\wZ, -)$, applied to the internal hom, gives back the set of morphisms
%
    \begin{equation}
   B_1(A)  =  \wAb(\wZ, A),    \qq    B_1(\Hom(B, C))  =  \wAb(B, C).
    \label{2.2.6}\end{equation}

	The unit-ball functor $ B_1 \c \wAb \to \Set $ has a left adjoint, associating to a set $ S $ the {\it free weighted abelian group} $ \wZ S $ generated by $ S$, namely the free abelian group generated by $ S $ with the weight
%
    \begin{equation}
	|\sum _x k_x.x|  =  \sum _x w(k_x),
    \label{2.2.7}\end{equation}
(where $ (k_x)_{x \in S}$ is a quasi-null family of integers). It is, of course, a sum of copies of $ \wZ$, the free weighted abelian group on one element, indexed on the elements of $ S$.

	Note that, in (\ref{2.2.3}), one {\it can} have $ |a \te  b| < |a|.|b| $ (as an obvious consequence of 
$ a\te b = (- a) \te  (- b)$, since we can have $ |a| < |- a|$.)

\subsection{Symmetry}\label{2.3}
The opposite weighted abelian group $ A\op$ has the same algebraic structure, with the {\it opposite weight}
%
    \begin{equation}
	|a|\op  =  |- a|    \qq    (a \in A).
    \label{2.3.1}\end{equation}

	A weighted abelian group is {\it symmetric} if it coincides with the opposite one, i.e. we always have $ |-a| = |a|$. Such objects form a full reflective and coreflective subcategory  $ !\wAb$.  The reflector  $ ! \c \wAb \to \,  !\wAb $ gives the symmetrised weighted abelian group $ !A$,  with the greatest symmetric weight $ || - || \le | - | $
%
    \begin{equation}
||a||  =  \inf_{\bf a} \; (\sum (|a_i| \mt |- a_i|),    \qq    ({\bf a} = (a_1,..., a_p), \; a_1+...+a_p= a).
    \label{2.3.4}\end{equation}

	The free symmetric weighted abelian group on one element is thus $ !\wZ$, with $ ||k|| $ the ordinary absolute value. It is easy to see that the reflector does not preserve finite products up to isometry, but only up to Lipschitz equivalence.	In \cite{G4}, we only used symmetric weighted abelian groups, called normed abelian groups; direction was obtained by enriching such object with a preorder, while here we prefer to derive the preorder from a (possibly non-symmetric) weight (\ref{2.4}).
	
	Formally, the notion of a symmetric weighted abelian group might seem to be preferable, since it amounts to an abelian group object in $  \wpSet$, while a weighted abelian group is just an abelian {\it monoid} in $  \wpSet $ which happens to be a group. But note that such problems arise whenever some form of `direction' is present; similarly, an ordered group is {\it not} a group in the category of ordered sets (since inversion must reverse the order), but an ordered monoid which happens to be a group.

\subsection{Other adjoints}\label{2.4}
The forgetful functor $ p \c \wAb \to \dAb $ with values in the category of preordered abelian groups equips a weighted group with the positive cone formed of the attainable elements (of finite weight)
%
    \begin{equation}
	a  \le  b    \q   \mbox{if}    \q    w(b - a)  <  \iy.
    \label{2.4.1}\end{equation}

	(Plainly, the functor $ p $ comes from the functor 
$ P \c \W \to \bt $ sending all $  \la < \iy $ to 1, see \ref{1.3.2}.) Its right adjoint $ wA $ gives to a preordered abelian group $ A $ the weight sending the positive cone to 0 and its complement to $\iy$.
	
	Enriching a previous result (\ref{2.2.7}), the forgetful functor $ \wAb \to \wSet $ has a left adjoint, associating to a weighted set $ X $ the {\it free weighted abelian group} $ \wZ X$, which is the free abelian group generated by the underlying set, equipped with the obvious weight
%
    \begin{equation}
|\sum _x k_x.x|  =  \sum _x w(k_x).|x|.
    \label{2.4.2}\end{equation}

\subsection{Weighted rings and modules}\label{2.5}
A {\it weighted ring} (with unit) will be a monoid in the monoidal category $ \wAb $ (\ref{2.2}). Since weighted abelian groups have an additive weight, and a multiplicatively-weighted tensor product, a weighted ring amounts to a ring $ R $ equipped with a weight function $ |a| \in \zety $ satisfying the following axioms (for all $ a, b \in R$)
%
    \begin{equation}
 |0|  =  0,   \q   |a + b|  \, \le \,  |a| + |b|,   \q
|1|  \le  1,  \q   |a.b|   \, \le \,   |a|.b|.
    \label{2.5.1}\end{equation}

	Note that the last condition should actually be written as $ |a.b| \le |a \te b|$, with respect to the norm of the tensor product $ R \te R $ of the underlying weighted abelian groups (cf. (\ref{2.2.3})). But, in the presence of the other axioms, these two conditions are equivalent (under universal quantifiers for $ a, b$, of course).

	In fact, the first condition is a trivial consequence of the second, since we always have $ |a \te b| \le |a|.b|$. Conversely, suppose we have $ |ab| \le |a|.|b| $ for all $ a, b$; if $ a \te b $ can be written as $ \sum a_i \te b_i$, it follows that $ ab = \sum  a_ib_i $ and
%
    \begin{equation}
\sum  |a_i|.|b_i|  \; \ge \;   \sum  |a_i b_i|  \; \ge \;   | \sum  a_i b_i |  \; = \;   |ab|,
    \label{2.5.2}\end{equation}
so that the greatest lower bound of such expressions $ \sum  |a_i|.|b_i|$ gives $ |a \te b| \ge |a.b|$.
	
	Typical examples are the weighted rings $ \wZ, \wQ, \wR $ of integers, rationals and reals, with
%
    \begin{equation}
   w(a)  =  a,   \for  a \ge 0,   \qq   w(a)  =  \iy,  \for  a < 0;
    \label{2.5.3}\end{equation}
their symmetrised versions $ !\wZ, !\wQ, !\wR $ are weighted by the ordinary absolute value.
	
	The category $ \wRng $ of weighted rings and contracting homomorphisms has thus forgetful functors with values in the categories $ \wAb $ and $ \wbMon $ of (additively) weighted abelian groups and multiplicatively weighted monoids.
	
	A {\it weighted module} $ A $ on the weighted ring $ R $ is a module with a weight satisfying
%
    \begin{equation}
 |0|  =  0,    \q    |a + b|  \le  |a| + |b|,    \q    | \la a|  \le  | \la|.|a|    \q    (a, b \in A;  \,   \la \in R).
    \label{2.5.4}\end{equation}

	Similarly for weighted vector spaces on weighted fields. Thus, a weighted abelian group is `the same' as a weighted module on the weighted ring $ \wZ$, while a symmetric weighted abelian group (\ref{2.3}) amounts to a weighted module on $ !\wZ$.
	
\section{Weighted categories}\label{3}
	Like monoids (\ref{2.1}), categories can be equipped with an additive weight or with a multiplicative one. In a weighted additive category (\ref{3.4}), defined as a category enriched on weighted abelian groups, the sum of parallel maps has a (sub)additive weight, while the weight of a composition is (sub)multiplicative.

\subsection{Categories with an additive weight}\label{3.1}
	An {\it additively weighted category} (called a normed category in \cite{Lw, BG}) will be a category $ X $ enriched on the symmetric monoidal closed category $  \wpSet$.
	
	Equivalently, $ X $ is a category and every morphism $ a $ is given a {\it weight}, or {\it cost} $ |a| \in \zety$, so that two obvious axioms are satisfied, for identities and composition:
%
    \begin{equation*}
    \begin{array}{lll}
 \mbox{(\wp cat.0)}  \;   &  |1_x|  =  0, &  \mbox{ for all objects }  x  \mbox{ of }  X,
\\[2pt]
 \mbox{(\wp cat.1)}    &  |ba|  \le  |a| + |b|,  &  \mbox{ for all pairs of consecutive arrows }  a, b. \qq\qq
    \label{3.1.1}    \end{array}
     \end{equation*}

	Loosely speaking, we can think of $ |a| $ as a cost, in some sense (length, duration, price, energy, etc.) which can `at worst' be added in a composition (typically, a concatenation of procedures). A weighted group - abelian or not - can be viewed as a weighted groupoid on one object (cf. \ref{2.2}).
	
	An additively weighted category will also be called a {\it \wp-category}. A {\it \wp-functor} $ f \c X \to Y $ is a functor between such categories, satisfying the condition $ |f(a)| \le |a|$, for all morphisms $ a $ of $ X$; a {\it \wp-transformation} $ \ph \c f \to g $ is a natural transformation between \wp-functors. All this forms the 2-category $  \wpCat $ of (small) \wp-categories.

	Some examples can be found in \cite{Lw}, p. 139. In a sequel, we will construct the {\it fundamental weighted category} of a generalised metric space, in the sense of Lawvere \cite{Lw}, or more generally of a `weighted space' (a topological space equipped with a weight function for continuous paths, under suitable axioms). Other examples are in \ref{3.3}.
	
	Incidentally, also the notion of a {\it sup-weighted category}, enriched on the {\it cartesian} monoidal structure $ (\wSet,  \ti )$, can be of interest. Then, the second axiom above has to be replaced with: $ |ba| \le |a|  \jn  |b|$. The previous `additive' notion is to the latter what a metric space is to an ultrametric one. 

\subsection{Categories with a multiplicative weight}\label{3.2}
Similarly, a {\it multiplicatively weighted category} will be a category $ X $ enriched over the symmetric monoidal closed category $  \wbSet$.

	Equivalently, $ X $ is a category and every morphism $ a $ is equipped with a {\it weight} $ |a| \in \zety$, so that:
%
    \begin{equation*}
    \begin{array}{lll}
 \mbox{(\wb cat.0)} \;   &  |1_x|  =  1, &  \mbox{ for all objects }  x  \mbox{ of }  X,
\\[2pt]
 \mbox{(\wb cat.1)}   &  |ba|  \le  |a|.|b|,  &  \mbox{ for all pairs of consecutive arrows }  a, b. \qq\qq
    \label{3.2.1}    \end{array}
     \end{equation*}

	We also speak of a {\it \wb-category}. A {\it \wb-functor} $ f \c X \to Y $ is a functor between such categories, satisfying the condition $ |f(a)| \le |a|$, for all morphisms $ a $ of $ X$; a {\it \wb-transformation} $ \ph \c f \to g $ is a natural transformation between \wb-functors. We have now the 2-category $ \wbCat $ of (small) \wb-categories.
	
	The categories of normed vector spaces and Banach spaces, with all linear maps (or with the bounded ones, or with linear contractions), have a classical multiplicative weight, the Lipschitz norm. Note that the norm of the identity of the null space is 0; this happens less exceptionally with {\it semi}normed vector space or weighted vector spaces. Other examples are considered below.

\subsection{Other examples}\label{3.3}
	The category of w-sets and {\it all} mappings has an {\it additive} weight $ |f|_0 $ as defined above for a mapping $ f \c X \to Y $ (see (\ref{1.5.3}))
%
    \begin{equation}
    \begin{array}{lll}
|f|_0  =  \sup_x (|f(x)| - |x|)  =  min  \{ \la \in \zety \;  |  \mbox{ for all }  x \in X,  |f(x)| \le  \la + |x| \},
\\[2pt]
|f|_0 \le  \la \;  \mbox{ if and only if, for all }   x \in X \c  \;   |f(x)| \le  \la + |x|,
    \label{3.3.1}    \end{array}
     \end{equation}
which identifies w-maps by the condition $ |f|_0 = 0$. It is thus trivial on $ \wSet$.

	But it has also a more usual {\it multiplicative} or {\it Lipschitz} weight (\ref{1.6.3})
%
    \begin{equation}
    \begin{array}{lll}
|f|_1  =  \sup_x (|f(x)| / |x|)  =  min  \{ \la \in \zety \;  |  \mbox{ for all }  x \in X,  |f(x)| \le  \la . |x| \},
\\[2pt]
|f|_1 \le  \la \;  \mbox{ if and only if, for all }   x \in X \c  \;   |f(x)| \le  \la . |x|,
    \label{3.3.2}    \end{array}
     \end{equation}
which identifies w-maps by the condition $ |f|_1 \le 1$, and is also of interest for $ \wSet$. There is a similar weight on the category of weighted abelian groups and {\it algebraic homomorphisms}; $ \wAb $ inherits a multiplicative weight.

\subsection{Weighted additive categories. }\label{3.4}
A {\it weighted additive category} $ A $ will be a category enriched on the monoidal category $ \wAb $ (\ref{2.2}). Extending the case of a weighted ring considered in \ref{2.5} (a weighted additive category on one object), all hom-sets $ A(X, Y) $ are weighted abelian groups, and composition is (sub)multiplicative
%
    \begin{equation}
    \begin{array}{lll}
|0|  =  0,   \q &   |f  + g|  \le  |f| + |g|  \qq &   (f, g \c X \to Y),
\\[2pt]
|1_X|  \le  1,   \q &   |gf| \le |f|.g|   \qq &   (f \c X \to Y,  \; g \c Y \to Z).
    \label{3.4.1}    \end{array}
     \end{equation}

	Also here the last condition is equivalent to $ |gf| \le |f \te g|$, with respect to the norm of the tensor product $ A(X, Y)  \te  A(Y, Z) $ (as in \ref{2.5}).

\subsection{Multiplicative norms and fibered categories}\label{3.5}
	We end with noting that Lawvere's article \cite{Lw} deals with the multiplicative aspect of Lipschitz norms in a different, much more general way.
	
	As hinted at there (see p. 150-151), one can form a category of generalised metric spaces, say $ {\bf M}$, fibered on the monoid of monoidal endofunctors of $ \Wp $ (as a category on one object). A morphism $ (f,  \la) \c X \to Y $ of $ {\bf M} $ consists of an arbitrary mapping $ f \c |X| \to |Y| $ between the underlying sets, together with a monoidal functor $  \la \c \Wp \to \Wp $ such that
%
    \begin{equation}
 \la (d_X (x, x'))  \; \ge \;   d_Y (f(x), f(x'))	\qq\qq		(x, x'  \in X),
    \label{3.5.1} \end{equation}
and the composition is obvious: $ (g,  \mu) \ci (f,  \la) = (gf,  \mu \la)$.

	Now, a monoidal functor $  \la \c \Wp \to \Wp $ is an increasing function such that:
%
    \begin{equation}
 \la \c \zety \to \zety,    \q    0  \ge   \la(0),    \q    \la(s) +  \la(t)  \ge   \la(s + t).
     \label{3.5.2} \end{equation}

	Every element $ \la \in \,  ]0, \iy[ $ gives a `linear endofunctor' 
	$ \hat{\la} (s) =  \la.s$, which, in condition (\ref{3.5.1}), says that $  \la $ is a Lipschitz constant for $ f$; composition multiplies such constants. But there are monoidal endofunctors, of interest in the theory of metric spaces, which are not `linear', e.g. the square root. Even restricting to `linear endofunctors', this approach is slightly different from the one we are following here; for instance, there are two candidates for $ \hat{0} $ (sending $\iy$ either to 0 or to $\iy$), but both give $\hat{\iy} \ci \hat{0} = \hat{0}$. 

	While this `fibered approach' is certainly important, the one we are following seems to give more directly the tools we need to develop `weighted' homology and homotopy, in a sequel.

\end{document}